\documentclass[12pt,draft]{article}
\usepackage{amsmath,amscd}
\usepackage{amssymb,latexsym,amsthm}
\usepackage{color}
\usepackage[spanish,english]{babel}
\usepackage{amssymb}
\usepackage{color}
\usepackage{amsmath,amsthm,amscd}
\usepackage[latin1]{inputenc}
\usepackage{lscape}
\usepackage{fancyhdr}
\usepackage{amsfonts}

\numberwithin{equation}{section}
\newtheorem{theorem}{Theorem}[section]

\newtheorem{proposition}[theorem]{Proposition}

\theoremstyle{definition}

\newtheorem{definition}[theorem]{Definition}
\newtheorem{examples}[theorem]{Examples}

\newtheorem{remark}[theorem]{Remark}

\newcommand{\cA}{\mbox{${\cal A}$}}

\newcommand{\cN}{\mbox{${\cal N}$}}

\newcommand{\cU}{\mbox{${\cal U}$}}

\newcommand{\cW}{\mbox{${\cal W}$}}

\title{\textbf{PBW bases for some 3-dimensional\\ skew polynomial algebras}}
\author{Armando Reyes\footnote{Departamento de Matem\'aticas. e-mail: mareyesv@unal.edu.co}\\ Universidad Nacional de
Colombia, Bogot\'a \\ H\'ector Su\'arez\footnote{Escuela de Matem\'aticas y Estad\'istica. e-mail: hector.suarez@uptc.edu.co}\\ Universidad Pedag\'ogica y Tecnol\'ogica de Colombia, Tunja}
\date{}
\begin{document}
\maketitle
\begin{abstract}
\noindent The aim of this paper is to establish necessary and sufficient algorithmic conditions to guarantee that an algebra is actually a 3-dimensional skew polynomial algebra in the sense of Bell and Smith \cite{SmithBell}.

\bigskip

\noindent \textit{Key words and phrases:} skew polynomial algebra, diamond lemma, skew PBW extension.

\bigskip

\noindent 2010 \textit{Mathematics Subject Classification.} 16S36, 16S32, 16S30.
\bigskip

\end{abstract}

\section{Introduction}
In the study of commutative and non-commutative algebras, it is important to specify one PBW (Poincar\'e-Birkhoff-Witt) basis for every one of them, since this allows us to characterize several properties with physical and mathematical meaning. This fact can be appreciated in several works. For instance, PBW theorem for the universal enveloping algebra of  a Lie algebra \cite{Dixmier1996}; PBW theorem for quantized universal enveloping algebras \cite{Yamane1989}; quantum PBW theorem for a wide class of associative algebras \cite{Berger1992}; PBW bases for quantum groups using the notion of Hopf algebra \cite{Ringel1996}, and others. With all these results in mind, in this article we wish to investigate a criteria and some algorithms which decide whether a given ring with some variables and relations can be expressed as a {\em 3-dimensional skew polynomial algebra} defined by Bell and Smith \cite{SmithBell} (Definition \ref{3dimensionaldimension}). We follow the original ideas by Bergman in \cite{Bergman1978} and the treatments established by Bueso et. al., \cite{BuesoGT2003} and Reyes \cite{ReyesPhD}. \\

The paper is organized as follows. Section \ref{Some3dimensionall} contains the algebras of interest for us in this paper, the 3-dimensional skew polynomial algebras. We recall its definition (Definition \ref{3dimensionaldimension}) and its classification (Proposition \ref{3-dimensionalClassification}). Section \ref{sectionDiamondLemma} treats  the definitions and preliminary results with the aim of establishing the important result of this paper (Theorem \ref{GomezTorrecillasTheorem3.21}). In Section \ref{SkewPoincareBirkhoffWittTheorem} we establish the algorithms which allow us to decide whether an algebraic structure, defined by variables and relations between them, can be considered as a 3-dimensional skew polynomial algebra (expressions (\ref{gordito}) - (\ref{pss6})). Finally, in Section  \ref{examples} we present some examples which illustrate the results obtained in Section \ref{sectionDiamondLemma}  and the algorithms formulated in Section \ref{SkewPoincareBirkhoffWittTheorem}.\\

Throughout this paper the letter $\Bbbk$ will denote a field.
\section{3-dimensional skew polynomial algebras}\label{Some3dimensionall}
The universal enveloping algebra $\cU(\mathfrak{sl}(2,\Bbbk))$ of
the  Lie algebra $\mathfrak{sl}(2,\Bbbk)$, the dispin algebra
$\cU(osp(1,2))$ and Woronowicz's algebra
$\cW_{\nu}(\mathfrak{sl}(2,\Bbbk))$ (see Examples \ref{chapeco}) are examples of algebras
classified by Bell and Smith in \cite{SmithBell}, which are known
as \textit{3-dimensional skew polynomial algebras}. These algebras are particular examples of a more general family of non-commutative rings known as {\em skew PBW extensions} or $\sigma$-{\em PBW extensions}. For these extensions several pro\-per\-ties have been characterized (for example, Noetherianess, regularity, Serre's Theorem, global homological, Krull, Goldie and Gelfand-Kirillov dimensions, Auslander's regularity, prime ideals, incomparability and prime length of prime ideals, higher algebraic $K$-theory, cyclic homology, Armendariz, Baer, quasi-Baer, p.p. and p.q.-Baer, and Koszul pro\-per\-ties, and other ring and module theoretical pro\-per\-ties, c.f. \cite{LezamaReyes2014}, \cite{LezamaAcostaReyes2015},  \cite{ReyesPhD}, \cite{Reyes2013}, \cite{Reyes2014}, \cite{Reyes2014UIS}, \cite{Reyes2015}, \cite{Reyes2018}, \cite{ReyesSuarez2016a}, \cite{ReyesSuarez2016b},  \cite{ReyesSuarez2016c}, \cite{ReyesSuarezClifford2017}, \cite{ReyesSuarezskewCY2017}, \cite{ReyesSuarezUMA2018}, \cite{ReyesYesica}, \cite{SuarezReyes2016}, \cite{SuarezReyesgenerKoszul2017} and others), which means that all these properties have been also investigated for 3-dimensional skew polynomial algebras. Nevertheless, since by definition (see Definition \ref{3dimensionaldimension} below) these algebras are required to have a PBW basis, we consider important to establish necessary and sufficient algorithmic conditions to guarantee that an algebra defined by generators and relations is precisely one of these skew polynomial algebras (this is done in Sections \ref{sectionDiamondLemma} and \ref{SkewPoincareBirkhoffWittTheorem}), and then can apply all results above. With this objective, we start recalling their definition and their characterization.
\begin{definition}[\cite{SmithBell}; \cite{Rosenberg1995}, Definition C4.3]\label{3dimensionaldimension}
\textit{A 3-dimensional skew polynomial algebra $\cA$} is a
$\Bbbk$-algebra generated by the variables $x,y,z$ restricted to relations $
yz-\alpha zy=\lambda,\ zx-\beta xz=\mu$, and $xy-\gamma
yx=\nu$, such that
\begin{enumerate}
\item $\lambda, \mu, \nu\in \Bbbk+\Bbbk x+\Bbbk y+\Bbbk z$, and $\alpha, \beta, \gamma \in \Bbbk^{*}$;
\item Standard monomials $\{x^iy^jz^l\mid i,j,l\ge 0\}$ are a $\Bbbk$-basis of the algebra.
\end{enumerate}
\end{definition}
\begin{remark}\label{Gererere}
If we consider the variables $x_1:=x,\ x_2:= y,\ x_3:=z$, then the relations established in Definition \ref{3dimensionaldimension} can be formulated in the following way:
\begin{align*}
x_3x_2 - \alpha^{-1} x_2x_3 = &\ r_0^{(2,3)} + r_1^{(2,3)}x_1 + r_2^{(2,3)}x_2 + r_3^{(2,3)}x_3,\\
x_3x_1 - \beta x_1x_3 = &\ r_0^{(1,3)} + r_1^{(1,3)}x_1 + r_2^{(1,3)}x_2 + r_3^{(1,3)}x_3,\\
x_2x_1 - \gamma^{-1}x_1x_2 = &\ r_0^{(1,2)} + r_1^{(1,2)}x_1 +
r_2^{(1,2)}x_2 + r_3^{(1,2)}x_3,
\end{align*}
where the elements $r's$ belong to the field $\Bbbk$.
\end{remark}
Next proposition establishes a classification of 3-dimensional
skew polynomial algebras.
\begin{proposition}[\cite{Rosenberg1995}, Theorem C.4.3.1, 2.5 in \cite{SmithBell}]\label{3-dimensionalClassification}
If $\cA$ is a 3-dimensional skew polynomial algebra, then $\cA$
is one of the following algebras:
\begin{enumerate}
\item [\rm (a)] if $|\{\alpha, \beta, \gamma\}|=3$, then $\cA$ is defined by the relations $yz-\alpha zy=0,\ zx-\beta xz=0,\ xy-\gamma yx=0$.
\item [\rm (b)] if $|\{\alpha, \beta, \gamma\}|=2$ and $\beta\neq \alpha =\gamma =1$, then $\cA$ is one of the following algebras:
\begin{enumerate}
\item [\rm (i)] $yz-zy=z,\ \ \ zx-\beta xz=y,\ \ \ xy-yx=x${\rm ;}
\item [\rm (ii)] $yz-zy=z,\ \ \ zx-\beta xz=b,\ \ \ xy-yx=x${\rm ;}
\item [\rm (iii)] $yz-zy=0,\ \ \ zx-\beta xz=y,\ \ \ xy-yx=0${\rm ;}
\item [\rm (iv)] $yz-zy=0,\ \ \ zx-\beta xz=b,\ \ \ xy-yx=0${\rm ;}
\item [\rm (v)] $yz-zy=az,\ \ \ zx-\beta xz=0,\ \ \ xy-yx=x${\rm ;}
\item [\rm (vi)] $yz-zy=z,\ \ \ zx-\beta xz=0,\ \ \ xy-yx=0$,
\end{enumerate}
where $a, b$ are any elements of $\Bbbk$. All nonzero values of $b$
give isomorphic algebras.
\item [\rm (c)] If $|\{\alpha, \beta, \gamma\}|=2$ and $\beta\neq \alpha=\gamma\neq 1$, then $\cA$ is one of the following algebras:
\begin{enumerate}
\item [\rm (i)] $yz-\alpha zy=0,\ \ \ zx-\beta xz=y+b,\ \ \ xy-\alpha yx=0${\rm ;}
\item [\rm (ii)] $yz-\alpha zy=0,\ \ \ zx-\beta xz=b,\ \ \ xy-\alpha yx=0$.
\end{enumerate}
In this case, $b$ is an arbitrary element of $\Bbbk$. Again, any
nonzero values of $b$ give isomorphic algebras.
\item [\rm (d)] If $\alpha=\beta=\gamma\neq 1$, then $\cA$ is the algebra defined by the relations  $yz-\alpha zy=a_1x+b_1,\ zx-\alpha xz=a_2y+b_2,\ xy-\alpha yx=a_3z+b_3$. If $a_i=0\ (i=1,2,3)$, then all nonzero values of $b_i$ give isomorphic
algebras.
\item [\rm (e)] If $\alpha=\beta=\gamma=1$, then $\cA$ is isomorphic to one of the following algebras:
\begin{enumerate}
\item [\rm (i)] $yz-zy=x,\ \ \ zx-xz=y,\ \ \ xy-yx=z${\rm ;}
\item [\rm (ii)] $yz-zy=0,\ \ \ zx-xz=0,\ \ \ xy-yx=z${\rm ;}
\item [\rm (iii)] $yz-zy=0,\ \ \ zx-xz=0,\ \ \ xy-yx=b${\rm ;}
\item [\rm (iv)] $yz-zy=-y,\ \ \ zx-xz=x+y,\ \ \ xy-yx=0${\rm ;}
\item [\rm (v)] $yz-zy=az,\ \ \ zx-xz=z,\ \ \ xy-yx=0${\rm ;}
\end{enumerate}
Parameters $a,b\in \Bbbk$ are arbitrary,  and all nonzero values of
$b$ generate isomorphic algebras.
\end{enumerate}
\end{proposition}
\section{Diamond lemma and PBW bases}\label{sectionDiamondLemma}
Bergman's Diamond Lemma  \cite{Bergman1978} provides a general method to prove that certain sets are bases of algebras which
are defined in terms of generators and relations. For instance,
the Poincar\'e-Birkhoff-Witt theorem, which appeared at first for
universal enveloping algebras of finite dimensional Lie algebras
(see \cite{Dixmier1996} for a detailed treatment) can be derived from it. PBW theorems have been considered several classes of commutative and noncommutative algebras (see \cite{Yamane1989}, \cite{Berger1992}, \cite{Ringel1996}, and others). With this in mind, in this section we establish a criteria and some algorithms which decide whether a given ring with some variables
and relations can be expressed as a 3-dimensional skew polynomial algebra in the sense of Definition \ref{3dimensionaldimension}. We follow the original ideas presented by Bergman \cite{Bergman1978} and the treatments developed by Bueso et. al. \cite{BuesoGT2003} and Reyes \cite{ReyesPhD}.
\begin{definition}\label{definitionBergmanmm}
\begin{enumerate}
\item [\rm (i)] {\rm Let} $X$ {\rm be a non-empty set and denote by} $
\langle X\rangle$ {\rm and} $\Bbbk \langle X\rangle$ {\rm the free monoid on} $X$ {\rm and the free associative} $\Bbbk$-{\rm ring on} $X$, {\rm respectively. A subset} $Q\subseteq \langle X\rangle\times \Bbbk \langle X\rangle$ {\rm is called a} \textit{reduction system} {\rm for} $\Bbbk\langle X\rangle$. {\rm An element} $\sigma=(W_{\sigma},f_{\sigma})\in Q$ {\rm has components} $W_{\sigma}$ {\rm a word in} $\langle X \rangle $ {\rm and} $f_{\sigma}$ {\rm a polynomial in}  $\Bbbk\langle X\rangle$. {\rm Note that every reduction system for} $\Bbbk\langle X\rangle$ {\rm defines a factor ring} $A=\Bbbk\langle X\rangle/I_Q$, {\rm with} $I_Q$ {\rm the two-sided ideal of} $\Bbbk\langle X\rangle$ {\rm generated by the polynomials}
$W_{\sigma}-f_{\sigma}$, {\rm with} $\sigma\in Q$.
\item [\rm (ii)] {\rm If} $\sigma$ {\rm is an element of a reduction system} $Q$ {\rm and} $A,B\in \langle X\rangle$, {\rm the} $\Bbbk$-{\rm linear endomorphism} $r_{A\sigma B}:\Bbbk\langle X\rangle \to \Bbbk\langle X\rangle$, {\rm which fixes all elements in the basis} $\langle X\rangle$ {\rm different
from} $AW_{\sigma}B$ {\rm and sends this particular element to}
$Af_{\sigma}B$ {\rm is called a} \textit{reduction} {\rm for} $Q$. {\rm If} $r$ {\rm is a reduction and} $f\in \Bbbk\langle X\rangle$, {\rm then} $f$ {\rm and}
$r(f)$ {\rm represent the same element in the} $\Bbbk$-{\rm ring} $\Bbbk\langle X\rangle/I_Q$. {\rm Thus, reductions may be viewed as rewriting rules
in this factor ring.}
\item [\rm (iii)] {\rm A reduction} $r_{A\sigma B}$ {\rm acts trivially on an element} $f\in \Bbbk\langle X\rangle$ {\rm if} $r_{A\sigma B}(f)=f$. {\rm An element} $f\in \Bbbk\langle X\rangle$ {\rm is said to be} irreducible {\rm under} $Q$ {\rm if all reductions act trivially on} $f$. {\rm Note that the set} $\Bbbk\langle X\rangle_{\rm irr}$ {\rm of all irreducible elements of} $\Bbbk\langle X\rangle$ {\rm under} $Q$ {\rm is a left submodule of} $\Bbbk\langle X\rangle$.
\item [\rm (iv)] {\rm Let} $f$ {\rm be an element of} $\Bbbk\langle X\rangle$. {\rm We say that} $f$ \textit{reduces} {\rm to} $g\in \Bbbk\langle X \rangle$, {\rm if there is a finite sequence} $r_1,\dotsc,r_n$ {\rm of reductions such that} $g=(r_n\dotsb r_1)(f)$. {\rm We will write} $f\to_Q g$. {\rm A finite sequence of reductions} $r_1,\dotsc,r_n$ {\rm is said to be}  final {\rm on} $f$, {\rm if} $(r_n\dotsb r_1)(f)\in \Bbbk\langle X\rangle_{\rm irr}$.
\item [\rm (v)] {\rm An element} $f\in \Bbbk\langle X\rangle$ {\rm is said to be} reduction-finite, {\rm if for every infinite sequence}  $r_1,r_2,\dotsc$ {\rm of reductions there exists some positive integer} $m$ {\rm such that} $r_i$ {\rm acts trivially on the element} $(r_{i-1}\dotsb r_1)(f)$, {\rm for every} $i>m$. {\rm If} $f$ {\rm is reduction-finite, then any maximal sequence of reductions} $r_1,\dotsc,r_n$ {\rm such that} $r_i$ {\rm acts non-trivially on the element}
$(r_{i-1}\dotsb r_1)(f)$, {\rm for} $1\le i\le n$, {\rm will be finite}. {\rm Thus, every reduction-finite element reduces to an irreducible element. We remark that
the set of all reduction-finite elements of} $\Bbbk\langle X\rangle$ {\rm is a left submodule of} $\Bbbk\langle X\rangle$.
\item [\rm (vi)] {\rm An element}  $f\in \Bbbk\langle X\rangle$ {\rm is said to be} \textit{reduction-unique} {\rm if it is reduction-finite and if its images under
all final sequences of reductions coincide. This value is denoted by} $r_Q(f)$.
\end{enumerate}
\end{definition}
\begin{proposition}[\cite{BuesoGT2003}, Lemma 3.13]\label{frooooooo}
{\rm (i)} The set $\Bbbk\langle X\rangle_{\rm un}$ of re\-duc\-tion-unique elements of $\Bbbk\langle X\rangle$ is a left submodule, and
$r_Q:\Bbbk\langle X\rangle_{\rm un}\to \Bbbk\langle X\rangle_{\rm irr}$
becomes an $\Bbbk$-linear map. {\rm (ii)} If $f,g,h\in \Bbbk\langle X\rangle$ are elements such that $ABC$ is reduction-unique for all terms $A,B,C$ occurring in respectively $f,g,h$, then $fgh$ is re\-duc\-tion-unique. Moreover, if $r$ is any reduction, then $fr(g)h$ is reduction-unique and $r_Q(fr(g)h)=r_Q(fgh)$.
\begin{proof}
(i) Consider $f, g \in \Bbbk\langle X\rangle_{\rm un},\ \lambda\in \Bbbk$. We know that $\lambda f + g$ is reduction-finite. Let $r_1,\dotsc, r_m$ be a sequence of reductions (note that it is final on this element), and $r:=r_m\dotsb r_1$ for the composition. Using that $f$ is reduction-unique, there is a finite composition of reductions $r'$ such that $(r'r)(f) = r_Q(f)$, and in a similar way, a composition of reductions $r''$ such that $(r''r'r)(g) = r_Q(g)$. Since $r(\lambda f + g)\in \Bbbk\langle X\rangle_{\rm irr}$, then $r(\lambda f + g) = (r''r'r)(\lambda f + g) = \lambda(r''r'r)(f) + (r''r'r)(g) = \lambda r_Q(f) + r_Q(g)$. Hence, the expression $r(\lambda f+g)$ is uniquely determined, and $\lambda f+g$ is reduction-unique. In fact, $r_Q(\lambda f+g) = \lambda r_Q(f) + r_Q(g)$, and therefore (i) is proved.

(ii) From (i) we know that $fgh$ is reduction-unique. Consider $r=r_{D\sigma E}$, for $\sigma\in Q,\ D, E\in \langle X\rangle$. The idea is to show that $fr(g)h$ is reduction-unique and $r_Q(fr(g)h) = r_Q(fgh)$. Note that if $f, g, h$ are terms $A, B, C$, then $r_{AD\sigma EC}(ABC) = Ar_{D\sigma E}(B)C$, that is, $Ar_{D\sigma E}(B)C$ is reduction-unique with the equality $r_Q(ABC) = r_Q(Ar_{D\sigma E}(B)C)$. Now, more generally, $f = \sum_{i}\lambda_i A_i,\ g=\sum_{j}\mu_j B_j,\ h=\sum_{k}\rho_kC_k$, where the indices $i, j, k$ run over finite sets, with  $\lambda_i, \mu_j, \rho_k$, and where $A_i, B_j, C_k$ are terms such that $A_iB_jC_k$ is reduction unique for every $i, j, k$. In this way, $fr(g)h = \sum_{i,j,k} \lambda_i\mu_j\rho_k A_ir(B_j)C_k$. Finally, since $ABC$ is reduction-finite,  for every $i, j, k$, and  $r_Q(A_ir(B_j)C_k) = r_Q(A_iB_jC_k)$, from (i), $fr(g)h$ is reduction-unique and $r_Q(fr(g)h) = r_Q(fgh)$.
\end{proof}
\end{proposition}
\begin{proposition}[\cite{BuesoGT2003}, Proposition 3.14]
If every element $f\in \Bbbk\langle X\rangle$ is re\-duc\-tion-finite
under a reduction system $Q$, and $I_Q$ is the ideal of $\Bbbk\langle
X\rangle$ generated by the set $\{W_{\sigma}-f_{\sigma}\mid \sigma
\in Q\}$ then $\Bbbk\langle X\rangle=\Bbbk\langle X\rangle_{irr}\oplus
I_Q$ if and only if every element of $\Bbbk\langle X\rangle$ is
reduction-unique.
\begin{proof}
Suppose that $\Bbbk\langle X\rangle = \Bbbk\langle X\rangle_{\rm irr}\oplus I_Q$ and consider $f\in \Bbbk\langle X\rangle$. Note that if $g, g'\in \Bbbk\langle X\rangle$ are elements for which $f$ reduces to $g$ and $g'$, then $g-g'\in \Bbbk\langle X\rangle\cap I_Q=\{0\}$, that is, $f$ is reduction-unique. Conversely, if every element of $\Bbbk\langle X\rangle$ is reduction-unique under $Q$, then $r_Q:\Bbbk\langle X\rangle \to \Bbbk\langle X\rangle_{\rm irr}$ is a $\Bbbk$-linear projection. Consider $f\in {\rm ker}(r_Q)$, that is, $r_Q(f)=0$. Then $f\in I_Q$, whence the ${\rm ker}(r_Q)\subseteq I_Q$, but in fact, ${\rm ker}(r_Q)$ contains $I_Q$: for every $\sigma\in Q, A, B\in \langle X\rangle$, we have $r_Q(A(W_{\sigma} - f_{\sigma})B) = r_Q(AW_{\sigma}B) - r_Q(Af_{\sigma}B) = 0$ from Proposition \ref{frooooooo}, when $r=r_{1\sigma 1}$.
\end{proof}
\end{proposition}
Under the previous assumptions, $A=\Bbbk\langle X\rangle/I_Q$ may be identified with the left free $\Bbbk$-module $\Bbbk\langle X\rangle_{irr}$ with $\Bbbk$-module
structure given by the multiplication $f*g=r_Q(fg)$.
\begin{definition}
{\rm An} overlap ambiguity {\rm for} $Q$ {\rm is a} $5$-{\rm tuple}
$(\sigma,\tau, A,B,C)$, {\rm where} $\sigma, \tau\in Q$ {\rm and} $A,B,C\in
\langle X\rangle\ \backslash\ \{1\}$ {\rm such that} $W_{\sigma}=AB$ {\rm and}
$W_{\tau}=BC$. {\rm This ambiguity is} solvable {\rm if there exist
compositions of reductions} $r,r'$ {\rm such that}
$r(f_{\sigma}C)=r'(Af_{\tau})$. {\rm Similarly,} {\rm a} 5-{\rm tuple} $(\sigma,\tau,A,B,C)$ {\rm with} $\sigma\neq \tau$ {\rm is called an}  inclusion ambiguity {\rm if} $W_{\tau}=B$ {\rm and} $W_{\sigma}=ABC$. {\rm This ambiguity is solvable if there are compositions of reductions}
$r,r'$ {\rm such that} $r(Af_{\tau}B)=r'(f_{\sigma})$.
\end{definition}
\begin{definition}
{\rm A partial monomial order} $\le$ {\rm on} $\langle X\rangle$ {\rm is said to be} compatible {\rm with} $Q$ {\rm if} $f_{\sigma}$ {\rm is a linear combination of
terms} $M$ {\rm with} $M<W_{\sigma}$, {\rm for all} $\sigma \in Q$.
\end{definition}
\begin{proposition}[\cite{BuesoGT2003}, Proposition 3.18]\label{GomezProposition3.18}
If $\le$ is a monomial partial order on $\langle X\rangle$
satisfying the descending chain condition and compatible with a
reduction system $Q$, then every element $f\in \Bbbk\langle X\rangle$
is reduction-finite. In particular, every element of $\Bbbk\langle
X\rangle$ reduces under $Q$ to an irreducible element.
\end{proposition}
Let $\le$ be a monoid partial order on $\langle X\rangle$
compatible with the reduction system $Q$. Let $M$ be a term in
$\langle X\rangle$ and write $Y_M$ for the submodule of $\Bbbk\langle
X\rangle$ spanned by all polynomials of the form
$A(W_{\sigma}-f_{\sigma})B$, where $A,B \in \langle X\rangle$ are
such that $AW_{\sigma}B<M$. We will denote by $V_M$ the submodule
of $\Bbbk\langle X \rangle$ spanned by all terms $M'< M$. Note that
$Y_M\subseteq V_M$.
\begin{definition}
{\rm An overlap ambiguity} $(\sigma,\tau, A,B,C)$ {\rm is said to be} resolvable {\rm relative to} $\le$ {\rm if} $f_{\sigma}C-Af_{\tau}\in
Y_{ABC}$. {\rm An inclusion ambiguity} $(\sigma,\tau,A,B,C)$ {\rm is said to
be} resolvable {\rm relative to} $\le$ {\rm if}
$Af_{\tau}C-f_{\sigma}\in Y_{ABC}$.
\end{definition}
If $r$ is a finite composition of reductions, and $f$ belongs to
$V_M$, then $f-r(f)\in Y_M$. Hence,  $f\in Y_M$ if and only if $r(f)\in Y_M$ (\cite{Reyes2013}, Proposition 3.1.8).\\

From the results above we obtain the important theorem of this section.
\begin{theorem}[Bergman's Diamond Lemma \cite{Bergman1978}; \cite{BuesoGT2003}, Theorem 3.21]\label{GomezTorrecillasTheorem3.21}
Let $Q$ be a reduction system for the free associative $\Bbbk$-ring
$\Bbbk\langle X\rangle$, and let $\le$ be a monomial partial order on
$\langle X\rangle$, compatible with $Q$ and satisfying the
descending chain condition. The following conditions are
equivalent: {\rm (i)} all ambiguities of $Q$ are resolvable; {\rm (ii)} all ambiguities of $Q$ are resolvable relative to $\le$; {\rm (iii)} all elements of $\Bbbk\langle X\rangle$ are reduction-unique under $Q$; {\rm (iv)}   $\Bbbk\langle X\rangle=\Bbbk\langle X\rangle_{\rm irr}\oplus I_Q$.
\end{theorem}
\section{Algorithms}\label{SkewPoincareBirkhoffWittTheorem}
Throughout this section we will consider the lexicographical degree order
$\preceq_{\rm deglex}$ to be defined on the variables $x_1,\dotsc,x_n$.
\begin{definition}
{\rm A reduction system} $Q$ {\rm for the free associative} $\Bbbk$-{\rm ring} given by $\Bbbk\langle
x_1,\dotsc,x_n\rangle$ {\rm is said to be a} $\preceq_{\rm
deglex}$-\textit{skew reduction system} {\rm if the following
conditions hold:} {\rm (i)} $Q=\{(W_{ji}, f_{ji})\mid 1\le i< j\le n\}$; {\rm (ii)} {\rm for every}  $j>i$, $W_{ji}=x_jx_i$ {\rm and}
$f_{ji}=c_{i,j}x_ix_j+p_{ji}$,
{\rm where} $c_{i,j}\in \Bbbk\ \backslash\ \{0\}$ {\rm and}  $p_{ji}\in \Bbbk\langle
x_1,\dotsc,x_n\rangle$; {\rm (iii)}  {\rm for each} $j>i$, ${\rm lm}(p_{ji}) \preceq_{\rm deglex} x_ix_j$. {\rm We will denote} $(Q,\preceq_{\rm deglex})$ {\rm this type of reduction systems.}
\end{definition}
Note that if $0\neq p\in
\sum_{\alpha}r_{\alpha}x^{\alpha}$, $r_{\alpha}\in \Bbbk$, we consider its
\textit{Newton diagram} as $ \cN(p):=\{\alpha\in \mathbb{N}^{n}\mid
r_{\alpha}\neq 0\}$.  Let ${\rm exp}(p):={\rm max}\ \cN(p)$. In
this way, by Proposition \ref{GomezProposition3.18} every element
$f\in \Bbbk\langle x_1,\dotsc,x_n\rangle$ reduces under $Q$ to an
irreducible element. Let $I_Q$ be the two-sided ideal of $\Bbbk\langle
x_1,\dotsc,x_n\rangle $ generated by $W_{ji}-f_{ji}$, for $1\le i<
j\le n$. If $x_i+I_Q$ is also represented by $x_i$, for each $1\le
i\le n$, then we call \textit{standard terms} in $A$. Proposition
\ref{GomezProposition4.3} below shows that any polynomial reduces
under $Q$ to some standard polynomial and hence standard terms in
$A$ generate this algebra as a left free $\Bbbk$-module.
\begin{proposition}[\cite{BuesoGT2003}, Lemma 4.2]\label{GomezLemma4.2}
If $(Q,\preceq_{\rm deglex})$ is a skew reduction system, then the set
$\Bbbk\langle x_1,\dotsc,x_n\rangle_{irr}$ is the left submodule of $\Bbbk\langle x_1,\dotsc,x_n\rangle$ consisting of all standard
polynomials $f\in \Bbbk\langle x_1,\dotsc,x_n\rangle$.
\begin{proof}
It is clear that every standard term is irreducible. Now, let us
see that if a monomial $M=\lambda x_{j_1}\dotsb x_{j_s}$ is not
standard, then some reduction will act non-trivially on it. If
$s<2$ the monomial is clearly standard. This is also true if
$j_k\le j_{k+1}$, for every $1\le k\le s-1$. Let $s\ge 2$. There
exists $k$ such that $j_k>j_{k+1}$ and $M=Cx_jx_iB=CW_{ji}B$ where
$j=j_k$, $i=j_{k+1}$ and where $C$ and $B$ are terms. Then $CW_{ji}B\to_Q Cf_{ji}B$ acts non trivially on $M$.
\end{proof}
\end{proposition}
\begin{proposition}[\cite{BuesoGT2003}, Proposition 4.3]\label{GomezProposition4.3}
If $(Q,\preceq_{\rm deglex})$ is a skew reduction system for $\Bbbk\langle x_1,\dotsc,x_n\rangle$, then every element of $\Bbbk\langle
x_1,\dotsc,x_n\rangle$ reduces under $Q$ to a standard polynomial.
Thus the standard terms in $A=\Bbbk\langle x_1,\dotsc,x_n\rangle/I_Q$
span $A$ as a left free module over $\Bbbk$.
\begin{proof}
It follows from Proposition \ref{GomezLemma4.2} and Proposition
\ref{GomezProposition3.18}.
\end{proof}
\end{proposition}
Next, we present an algorithm to reduce any polynomial in $\Bbbk\langle
x_1,\dotsc,x_n\rangle$ to its standard representation modulo
$I_Q$. The basic step in this algorithm is the reduction of terms
to polynomials of smaller leading term. In the proof of Proposition
\ref{GomezLemma4.2} we can choose $k$ to be the least integer such
that $j_k>j_{k+1}$, thus yielding a procedure to define for every
non-standard monomial $\lambda M$ a reduction denoted ${\rm red}$
that acts non-trivially on $M$. In this way, the linear map ${\rm
red}:\Bbbk\langle x_1,\dotsc, x_n\rangle \to \Bbbk\langle x_1,\dotsc,
x_n\rangle $ depends on $M$. However, the following procedure is
an algorithm.
{\scriptsize{
\begin{center}
\begin{tabular}{|p{10.5cm}|}\hline
\textsf{}\\

\centerline{\textbf{Algorithm: Monomial reduction algorithm}}\label{algorithm1}\\

\setlength{\parindent}{1cm}\textbf{INPUT:} $M=\lambda x_{j_1}\dotsb x_{j_r}$ a non standard monomial.\\

\setlength{\parindent}{1cm}\textbf{OUTPUT:} $p={\rm red}(M)$, a reduction under $Q$ of the monomial $M$\\

\setlength{\parindent}{1cm}\textbf{INITIALIZATION:} $k=1, C=\lambda$\\

\ \ \ \ \ \ \ \ \ \ \ \ \setlength{\parindent}{1cm}\textbf{WHILE} $j_k \le j_{k+1}$ \textbf{DO}\\

\ \ \ \ \ \ \ \ \ \ \ \ $C=Cx_{j_k}$\\

\ \ \ \ \ \ \ \ \ \ \ \ $k=k+1$\\

\ \ \ \ \ \ \ \ \ \ \ \ \ \ \textbf{IF} $k+2\le r$ \textbf{THEN}\\

\ \ \ \ \ \ \ \ \ \ \ \ \ \ \ \ \ $B=x_{j_{k+2}}\dotsb x_{j_r}$\\

\ \ \ \ \ \ \ \ \ \ \ \ \ \ \textbf{ELSE}\\

\ \ \ \ \ \ \ \ \ \ \ \ \ \ \ \ \ $B=1$\\

\ \ \ \ \ \ \ \ \ \ \ \ \ \ \ \ \ $j=j_k$, $i=j_{k+1}$\\

\ \ \ \ \ \ \ \ \ \ \ \ \ \ \ \ \ $p=Cf_{j,i}B$.
\vspace{0.5cm}
\\ \hline
\end{tabular}
\end{center}}}
An element $f\in \Bbbk\langle x_1,\dotsc,
x_n\rangle$ is called \textit{normal} if $ {\rm
deg}(X_t)\preceq_{\rm deglex} {\rm deg}({\rm lt}(f))$, for every
term $X_t\neq {\rm lt}(f)$ in $f$.
\begin{proposition}[\cite{BuesoGT2003}, Proposition 4.5]\label{GomezTorrecillasProposition4.5}
If $(Q,\preceq_{\rm deglex})$ is a skew reduction system, then there exists a $\Bbbk$-linear map
${\rm stred}_Q:\Bbbk\langle x_1,\dotsc, x_n\rangle \to \Bbbk\langle
x_1,\dotsc, x_n\rangle_{ \rm irr}$
satisfying the following conditions: {\rm (i)} for every element $f$ of $\Bbbk\langle x_1,\dotsc, x_n\rangle$, there exists a finite sequence $r_1,\dotsc$, $r_m$ of reductions such that ${\rm stred}_Q(f)=(r_m\dotsb r_1)(f)$; {\rm (ii)}  if $f$ is normal, then we obtain ${\rm deg}({\rm lm}(f))={\rm deg}({\rm lm}({\rm stred}_Q(f)))$.
\end{proposition}
From the proof of Proposition \ref{GomezTorrecillasProposition4.5}
we obtain the next algorithm. Remark \ref{valverde} and Theorem \ref{GomezTorrecillas2Theorem 4.7} are the key results connecting this section with 3-dimensional skew polynomial algebras. \\
{\scriptsize{\begin{center}
\begin{tabular}{|p{10.2cm}|}\hline
\textsf{}

\centerline{\textbf{Algorithm: Reduction to standard form algorithm}} \label{algorithm2}\\

\setlength{\parindent}{1cm}\textbf{INPUT:} $f$ a non-standard polynomial.\\

\setlength{\parindent}{1cm}\textbf{OUTPUT:} $g={\rm stred}_Q(f)$ a standard reduction under $Q$ of $f$\\

\setlength{\parindent}{1cm}\textbf{INITIALIZATION:} $g=0$\\

\ \ \ \ \ \ \ \ \ \ \ \setlength{\parindent}{1cm}\textbf{WHILE} $f\neq 0$ \textbf{DO}\\

\ \ \ \ \ \ \ \ \ \ \ \ \ \ \textbf{IF} ${\rm lm}(f)$ is standard \textbf{THEN}\\

\ \ \ \ \ \ \ \ \ \ \ \ \ \ \ \ \ $f=f-{\rm lm}(f)$\\

\ \ \ \ \ \ \ \ \ \ \ \ \ \ \ \ \ $g=g+{\rm lm}(g)$\\

\ \ \ \ \ \ \ \ \ \ \ \ \ \ \textbf{ELSE}\\

\ \ \ \ \ \ \ \ \ \ \ \ \ \ \ \ \ $f=f-{\rm lm}(f)+{\rm red}({\rm
lm}(f))$.
\vspace{0.5cm}
\\ \hline
\end{tabular}
\end{center}}}
\begin{remark}\label{valverde}
A free left $\Bbbk$-module $A$ is a 3-dimensional skew polynomial algebra with respect
to $\preceq_{\rm deglex}$ if and only if it is isomorphic to $\Bbbk\langle x_1,\dotsc,x_n\rangle /I_Q$, where $Q$ is a skew reduction system with respect to $\preceq_{\rm deglex}$.
\end{remark}
By Theorem \ref{GomezTorrecillasTheorem3.21}, the set of all standard terms forms a $\Bbbk$-basis for $A$ given by $A=\Bbbk\langle
x_1,\dotsc, x_n\rangle/I_Q$. We have the following key result:
\begin{theorem}[\cite{BuesoGT2003}, Theorem 4.7]\label{GomezTorrecillas2Theorem 4.7}
Let $(Q,\preceq_{\rm deglex})$ be a skew reduction system on
$\Bbbk\langle x_1,\dotsc, x_n\rangle$ and let $A=\Bbbk\langle
x_1,\dotsc,x_n\rangle/I_Q$. For $1\le i<j<k\le n$, let $g_{kji},
h_{kji}$ be elements in $\Bbbk\langle x_1,\dotsc, x_n\rangle$ such
that $x_kf_{ji}$ {\rm (}resp.  $f_{kj}x_i${\rm )} reduces to
$g_{kji}$ {\rm (}resp. $h_{kji}${\rm )} under $Q$. The following
conditions are equivalent:
\begin{enumerate}
\item [\rm (i)] $A$ is a 3-dimensional skew polynomial algebra over $\Bbbk$;
\item [\rm (ii)] the standard terms form a basis of $A$ as a left free $\Bbbk$-module;
\item [\rm (iii)] $g_{kji}=h_{kji}$,  for every $1\le i<j<k\le n$;
\item [\rm (iv)] ${\rm stred}_Q(x_kf_{ji})={\rm stred}_Q(f_{kj}x_i)$, for every $1\le i<j<k\le n$.
\end{enumerate}
Moreover, if $A$ is a 3-dimensional skew polynomial algebra, then ${\rm
stred}_Q=r_Q$ and $A$ is isomorphic as a left module to $\Bbbk\langle
x_1,\dotsc,x_n\rangle_{\rm irr}$ whose module structure is given
by the product $f*g:=r_Q(fg)$, for every $f,g\in \Bbbk\langle
x_1,\dotsc,x_n\rangle_{\rm irr}$.
\begin{proof}
The equivalence between (i) and (ii) as well between (i) and (iii)
is given by Theorem \ref{GomezTorrecillasTheorem3.21}. The
equivalence between (i) and (iv) is obtained from Theorem
\ref{GomezTorrecillasTheorem3.21} and Proposition
\ref{GomezTorrecillasProposition4.5}. The remaining statements are
also consequences of Theorem \ref{GomezTorrecillasTheorem3.21}.
\end{proof}
\end{theorem}
Theorem \ref{GomezTorrecillas2Theorem 4.7} gives an algorithm to
check whether $\Bbbk\langle x_1,\dotsc, x_n\rangle/I_Q$ is a skew PBW extension since ${\rm stred}_Q(x_kf_{ji})$ and ${\rm
stred}_Q(f_{kj}x_i)$ can be computed by means of Algorithm
\textquotedblleft Reduction to standard form
algorithm\textquotedblright.
\section{Examples}\label{examples}
Next, we consider Theorem \ref{GomezTorrecillas2Theorem 4.7} with the aim of showing the relations between the elements $r's$ which guarantee that one can have a 3-dimensional skew polynomial algebra with basis  given by Definition \ref{3dimensionaldimension}. If $x_1\prec x_2\prec x_3$ with the notation in Remark \ref{Gererere}, then $(Q,\preceq_{\rm deglex})$ is a skew reduction system and
{\small{\begin{align*}
{\rm stred}_Q(x_3f_{21}) = &\ x_3(\gamma^{-1} x_1x_2 + r_0^{(1,2)} + r_1^{(1,2)}x_1 + r_2^{(1,2)}x_2 + r_3^{(1,2)}x_3)\\
= &\ \gamma^{-1} x_3x_1x_2 + r_0^{(1,2)}x_3 + r_1^{(1,2)}x_3x_1 + r_2^{(1,2)}x_3x_2 + r_3^{(1,2)}x_3^{2}\\
= &\ \gamma^{-1} (\beta x_1x_3 + r_0^{(1,3)} + r_1^{(1,3)}x_1 + r_2^{(1,3)}x_2 + r_3^{(1,3)}x_3)x_2 + r_0^{(1,2)}x_3 \\
+ &\ r_1^{(1,2)}(\beta x_1x_3 + r_0^{(1,3)} + r_1^{(1,3)}x_1 + r_2^{(1,3)}x_2 + r_3^{(1,3)}x_3) \\
+ &\ r_2^{(1,2)}(\alpha^{-1}  x_2x_3 + r_0^{(2,3)} + r_1^{(2,3)}x_1 + r_2^{(2,3)}x_2 + r_3^{(2,3)}x_3) + r_3^{(1,3)}x_3^{2}\\
= &\ \gamma^{-1}  \beta x_1x_3x_2 + \gamma^{-1}  r_0^{(1,3)}x_2 + \gamma^{-1}  r_1^{(1,3)}x_1x_2 + \gamma^{-1}  r_2^{(1,3)}x_2^{2} + \gamma^{-1}  r_3^{(1,3)}x_3x_2 \\
+ &\ r_0^{(1,2)}x_3 + r_1^{(1,2)}\beta x_1x_3 + r_1^{(1,2)}r_0^{(1,3)} + r_1^{(1,2)}r_1^{(1,3)}x_1 + r_1^{(1,2)}r_2^{(1,3)}x_2 \\
+ &\ r_1^{(1,2)}r_3^{(1,3)}x_3 + r_2^{(1,2)}\alpha^{-1}  x_2x_3 + r_2^{(1,2)}r_0^{(2,3)} + r_2^{(1,2)}r_1^{(2,3)}x_1\\
+ &\ r_2^{(1,2)}r_2^{(2,3)}x_2 + r_2^{(1,2)}r_3^{(1,3)}x_3 + r_3^{(1,3)}x_3^{2}\\
= &\ \gamma^{-1}  \beta x_1(\alpha^{-1}  x_2x_3 + r_0^{(2,3)} + r_1^{(2,3)}x_1 + r_2^{(2,3)}x_2 + r_3^{(2,3)}x_3) + \gamma^{-1}  r_0^{(1,3)}x_2 \\
+ &\ \gamma^{-1}  r_1^{(1,3)}x_1x_2 + \gamma^{-1}  r_2^{(1,3)}x_2^{2} + \gamma^{-1}  r_3^{(1,3)}(\alpha^{-1}  x_2x_3 + r_0^{(2,3)} + r_1^{(2,3)}x_1 \\
+ &\ r_2^{(2,3)}x_2 + r_3^{(2,3)}x_3) + r_0^{(1,2)}x_3 + r_1^{(1,2)}\beta x_1x_3 +  r_1^{(1,2)}r_0^{(1,3)} \\
+ &\ r_1^{(1,2)}r_1^{(1,3)}x_1 + r_1^{(1,2)}r_2^{(1,3)}x_2 + r_1^{(1,2)}r_3^{(1,3)}x_3 + r_2^{(1,2)}\alpha^{-1}  x_2x_3 \\
+ &\ r_2^{(1,2)}r_0^{(2,3)} + r_2^{(1,2)}r_1^{(2,3)}x_1 + r_2^{(1,2)}r_2^{(2,3)}x_2 + r_2^{(1,2)}r_3^{(1,3)}x_3 + r_3^{(1,3)}x_3^{2}\\
= &\ \gamma^{-1}  \beta \alpha^{-1}  x_1x_2x_3 + \gamma^{-1}  \beta r_0^{(2,3)}x_1 +  \gamma^{-1}  \beta r_1^{(2,3)}x_1^{2} +  \gamma^{-1}  \beta r_2^{(2,3)}x_1x_2 \\
+ &\ \gamma^{-1}  \beta r_3^{(2,3)}x_1x_3 + \gamma^{-1}  r_0^{(1,3)}x_2 + \gamma^{-1}  r_1^{(1,3)}x_1x_2 \\
+ &\ \gamma^{-1}  r_2^{(1,3)}x_2^{2} + \gamma^{-1}  r_3^{(1,3)}\alpha^{-1}  x_2x_3 + \gamma^{-1}  r_3^{(1,3)}r_0^{(2,3)} + \gamma^{-1}  r_3^{(1,3)}r_1^{(2,3)}x_1 \\
+ &\ \gamma^{-1}  r_3^{(1,3)}r_2^{(2,3)}x_2 + \gamma^{-1}  r_3^{(1,3)} r_3^{(2,3)}x_3 + r_0^{(1,2)}x_3 + r_1^{(1,2)}\beta x_1x_3 \\
+ &\ r_1^{(1,2)}r_0^{(1,3)} + r_1^{(1,2)}r_1^{(1,3)}x_1 + r_1^{(1,2)}r_2^{(1,3)}x_2 + r_1^{(1,2)}r_3^{(1,3)}x_3\\
+ &\ r_2^{(1,2)}\alpha^{-1}  x_2x_3 + r_2^{(1,2)}r_0^{(2,3)} + r_2^{(1,2)}r_1^{(2,3)}x_1\\
+ &\ r_2^{(1,2)}r_2^{(2,3)}x_2 + r_2^{(1,2)}r_3^{(1,3)}x_3 +
r_3^{(1,3)}x_3^{2},
\end{align*}}}
or equivalently,
{\small{\begin{align*}
{\rm stred}_Q(x_3f_{21}) = &\ \gamma^{-1}  \beta \alpha^{-1}  x_1x_2x_3 + (\gamma^{-1}  \beta r_0^{(2,3)}  + \gamma^{-1}  r_3^{(1,3)}r_1^{(2,3)} + r_1^{(1,2)} r_1^{(1,3)}\\
+ &\ r_2^{(1,2)}r_1^{(2,3)})x_1 + (\gamma^{-1}  r_0^{(1,3)} + \gamma^{-1}  r_3^{(1,3)}r_2^{(2,3)} + r_1^{(1,2)}r_2^{(1,3)} + r_2^{(1,2)}r_2^{(2,3)})x_2 \\
+ &\ (\gamma^{-1}  r_3^{(1,3)}r_3^{(2,3)} + r_0^{(1,2)} + r_1^{(1,2)}r_3^{(1,3)} + r_2^{(1,2)}r_3^{(1,3)})x_3\\
+ &\ (\gamma^{-1}  r_1^{(1,3)} + \gamma^{-1}  \beta r_2^{(2,3)})x_1x_2 \\
+ &\ (\gamma^{-1}  \beta r_3^{(2,3)} + \beta r_1^{(1,2)})x_1x_3 + (\gamma^{-1}  \alpha^{-1}  r_3^{(1,3)} + \alpha^{-1}  r_2^{(1,2)})x_2x_3 \\
+ &\ \gamma^{-1}  \beta r_1^{(2,3)}x_1^{2} +  \gamma^{-1}  r_2^{(1,3)}x_2^{2} + r_3^{(1,3)}x_3^{2}\\
+ &\ \gamma^{-1}  r_3^{(1,3)}r_0^{(2,3)} + r_1^{(1,2)}r_0^{(1,3)} + r_2^{(1,2)}r_0^{(2,3)}.
\end{align*}}}
Next, we compute ${\rm stred}_Q(f_{32}x_1)$:
{\small{\begin{align*}
{\rm stred}_Q(f_{32}x_1) = &\ (\alpha^{-1}  x_2x_3 + r_0^{(2,3)} + r_1^{(2,3)}x_1 + r_2^{(2,3)}x_2 + r_3^{(2,3)}x_3)x_1\\
= &\ \alpha^{-1}  x_2x_3x_1 + r_0^{(2,3)}x_1 + r_1^{(2,3)}x_1^{2} + r_2^{(2,3)}x_2x_1 + r_3^{(2,3)}x_3x_1\\
= &\ \alpha^{-1}  x_2(\beta x_1x_3 + r_0^{(1,3)} + r_1^{(1,3)}x_1 + r_2^{(1,3)}x_2 + r_3^{(1,3)}x_3) + r_0^{(2,3)}x_1 \\
+ &\ r_1^{(2,3)}x_1^{2} + r_2^{(2,3)}(\gamma^{-1}  x_1x_2 + r_0^{(1,2)} + r_1^{(1,2)}x_1 + r_2^{(1,2)}x_2 + r_3^{(1,2)}x_3)\\
+ &\ r_3^{(2,3)}(\beta x_1x_3 + r_0^{(1,3)} + r_1^{(1,3)}x_1 + r_2^{(1,3)}x_2 + r_3^{(1,3)}x_3)\\
= &\ \beta \alpha^{-1}  x_2x_1x_3 + \alpha^{-1}  r_0^{(1,3)}x_2\\
+ &\ \alpha^{-1}  r_1^{(1,3)}x_2x_1 + \alpha^{-1}  r_2^{(1,3)}x_2^{2} + \alpha^{-1}  r_3^{(1,3)}x_2x_3 \\
+ &\ r_0^{(2,3)}x_1 + r_1^{(2,3)}x_1^{2} + \gamma^{-1}  r_2^{(2,3)}x_1x_2 + r_0^{(1,2)}r_2^{(2,3)}\\
+ &\ r_1^{(1,2)}r_2^{(2,3)}x_1 + r_2^{(1,2)}r_2^{(2,3)}x_2 + r_2^{(2,3)}r_3^{(1,2)}x_3 + \beta r_3^{(2,3)}x_1x_3 \\
+ &\ r_0^{(1,3)}r_3^{(2,3)} + r_1^{(1,3)}r_2^{(2,3)}x_1 + r_2^{(1,3)}r_3^{(2,3)}x_2 + r_3^{(1,3)}r_3^{(2,3)}x_3\\
= &\ \beta \alpha^{-1}  (\gamma^{-1}  x_1x_2 + r_0^{(1,2)} + r_1^{(1,2)}x_1 + r_2^{(1,2)}x_2 + r_3^{(1,2)}x_3)x_3 + \alpha^{-1}  r_0^{(1,3)}x_2 \\
+ &\ \alpha^{-1}  r_1^{(1,3)}(\gamma^{-1}  x_1x_2 + r_0^{(1,2)} + r_1^{(1,2)}x_1 + r_2^{(1,2)}x_2 + r_3^{(1,2)}x_3) + \alpha^{-1}  r_2^{(1,3)}x_2^{2}\\
+ &\ \alpha^{-1}  r_3^{(1,3)}x_2x_3 + r_0^{(2,3)}x_1 + r_1^{(2,3)}x_1^{2} + \gamma^{-1}  r_2^{(2,3)}x_1x_2 + r_0^{(1,2)}r_2^{(2,3)} \\
+ &\ r_1^{(1,2)}r_2^{(2,3)}x_1 + r_2^{(1,2)}r_2^{(2,3)}x_2 + r_2^{(2,3)}r_3^{(1,2)}x_3 + \beta r_3^{(2,3)}x_1x_3\\
+ &\ r_0^{(1,3)}r_3^{(2,3)} + r_1^{(1,3)}r_2^{(2,3)}x_1 + r_2^{(1,3)}r_3^{(2,3)}x_2 + r_2^{(1,3)}r_3^{(2,3)}x_3
\end{align*}}}
or what is the same,
{\small{\begin{align*}
{\rm stred}_Q(f_{32}x_1) = &\ \gamma^{-1}  \beta \alpha^{-1}  x_1x_2x_3 + \beta \alpha^{-1}  r_0^{(1,2)}x_3 + \beta \alpha^{-1}  r_1^{(1,2)}x_1x_3 + \beta \alpha^{-1}  r_2^{(1,2)}x_2x_3 \\
+ &\ \beta \alpha^{-1}  r_3^{(1,2)}x_3^{2} + \alpha^{-1}  r_0^{(1,3)}x_2 + \gamma^{-1}  \alpha^{-1}  r_1^{(1,3)}x_1x_2 + \alpha^{-1}  r_0^{(1,2)}r_1^{(1,3)} \\
+ &\ \alpha^{-1}  r_1^{(1,2)}r_1^{(1,3)}x_1 + \alpha^{-1}  r_1^{(1,3)}r_2^{(1,2)}x_2 + \alpha^{-1}  r_1^{(1,3)} r_3^{(1,2)}x_3 \\
+ &\ \alpha^{-1}  r_2^{(1,3)}x_2^{2} + \alpha^{-1}  r_3^{(1,3)}x_2x_3 + r_0^{(2,3)}x_1 + r_1^{(2,3)}x_1^{2} + \gamma^{-1}  r_2^{(2,3)}x_1x_2 \\
+ &\ r_0^{(1,2)}r_2^{(2,3)} + r_1^{(1,2)}r_2^{(2,3)}x_1 + r_2^{(1,2)}r_2^{(2,3)}x_2 +  r_2^{(2,3)}r_3^{(1,2)}x_3 \\
+ &\ \beta r_3^{(2,3)}x_1x_3 + r_0^{(1,3)}r_3^{(2,3)} + r_1^{(1,3)}r_2^{(2,3)}x_1 \\
+ &\ r_2^{(1,3)}r_3^{(2,3)}x_2 + r_2^{(1,3)}r_3^{(2,3)}x_3\\
= &\ \gamma^{-1}  \beta \alpha^{-1}  x_1x_2x_3 + (\alpha^{-1}  r_1^{(1,2)}r_1^{(1,3)} + r_0^{(2,3)} + r_1^{(1,2)}r_2^{(2,3)} + r_1^{(1,3)}r_2^{(2,3)})x_1\\
+ &\ (\alpha^{-1}  r_0^{(1,3)} + \alpha^{-1}  r_1^{(1,3)}r_2^{(1,2)} + r_2^{(1,2)}r_2^{(2,3)} + r_2^{(1,3)}r_3^{(2,3)})x_2\\
+ &\ (\beta \alpha^{-1}  r_0^{(1,2)} + \alpha^{-1}  r_1^{(1,3)}r_3^{(1,2)} + r_2^{(2,3)}r_3^{(1,2)} + r_2^{(1,3)}r_3^{(2,3)})x_3\\
+ &\ (\gamma^{-1}  \alpha^{-1}  r_1^{(1,3)} + \gamma^{-1}  r_2^{(2,3)})x_1x_2 + (\beta \alpha^{-1}  r_1^{(1,2)} + \beta r_3^{(2,3)})x_1x_3 \\
+ &\ (\beta \alpha^{-1}  r_2^{(1,2)} + \alpha^{-1}  r_3^{(1,3)})x_2x_3 + r_1^{(2,3)}x_1^{2} + \alpha^{-1}  r_2^{(1,3)}x_2^{2} + \beta \alpha^{-1}  r_3^{(1,2)}x_3^{2}\\
+ &\ \alpha^{-1}  r_0^{(1,2)}r_1^{(1,3)} + r_0^{(1,2)}r_2^{(2,3)} + r_0^{(1,3)}r_3^{(2,3)}.
\end{align*}}}
Since we need to satisy the relation ${\rm stred}_Q(x_3f_{21}) = {\rm stred}_Q(f_{32}x_1)$, the following equalities are necessary and sufficient to guarantee that an algebra generated by three variables (where coefficients commute with variables) can be considered as a 3-dimensional skew polynomial algebra in the sense of \cite{SmithBell}:
\begin{multline}\label{gordito}
\gamma^{-1}  \beta r_0^{(2,3)}  + \gamma^{-1}  r_3^{(1,3)}r_1^{(2,3)} + r_1^{(1,2)} r_1^{(1,3)} + r_2^{(1,2)}r_1^{(2,3)} = \\ \alpha^{-1}  r_1^{(1,2)}r_1^{(1,3)} + r_0^{(2,3)} + r_1^{(1,2)}r_2^{(2,3)} + r_1^{(1,3)}r_2^{(2,3)},
\end{multline}
\begin{multline}\label{flaquito}
\gamma^{-1}  r_0^{(1,3)} + \gamma^{-1}  r_3^{(1,3)}r_2^{(2,3)} + r_1^{(1,2)}r_2^{(1,3)} + r_2^{(1,2)}r_2^{(2,3)} = \\ \alpha^{-1}  r_0^{(1,3)} + \alpha^{-1}  r_1^{(1,3)}r_2^{(1,2)} + r_2^{(1,2)}r_2^{(2,3)} + r_2^{(1,3)}r_3^{(2,3)},
\end{multline}
\begin{multline}\label{moder}
\gamma^{-1}  r_3^{(1,3)}r_3^{(2,3)} + r_0^{(1,2)} + r_1^{(1,2)}r_3^{(1,3)} + r_2^{(1,2)}r_3^{(1,3)} = \\ \beta \alpha^{-1}  r_0^{(1,2)} + \alpha^{-1}  r_1^{(1,3)}r_3^{(1,2)} + r_2^{(2,3)}r_3^{(1,2)} + r_2^{(1,3)}r_3^{(2,3)},
\end{multline}
and
\begin{align}
\gamma^{-1}  r_1^{(1,3)} + \gamma^{-1}  \beta r_2^{(2,3)} = &\ \gamma^{-1}  \alpha^{-1}  r_1^{(1,3)} + \gamma^{-1}  r_2^{(2,3)} \label{pss1} \\
\gamma^{-1}  \beta r_3^{(2,3)} + \beta r_1^{(1,2)} = &\ \beta \alpha^{-1}  r_1^{(1,2)} + \beta r_3^{(2,3)} \label{pss2} \\
\gamma^{-1}  \alpha^{-1}  r_3^{(1,3)} + \alpha^{-1}  r_2^{(1,2)} = &\ \beta \alpha^{-1}  r_2^{(1,2)} + \alpha^{-1}  r_3^{(1,3)} \label{pss3} \\
\gamma^{-1}  \beta r_1^{(2,3)} = &\ r_1^{(2,3)}\label{cattt} \\
\gamma^{-1}  r_2^{(1,3)} = &\ \alpha^{-1}  r_2^{(1,3)} \label{doggg}\\
r_3^{(1,3)} = &\ \beta \alpha^{-1}  r_3^{(1,2)}\label{delfi} \\
\gamma^{-1}  r_3^{(1,3)}r_0^{(2,3)} + r_1^{(1,2)}r_0^{(1,3)} + r_2^{(1,2)}r_0^{(2,3)} = &\ \alpha^{-1}  r_0^{(1,2)}r_1^{(1,3)} + r_0^{(1,2)}r_2^{(2,3)} + r_0^{(1,3)}r_3^{(2,3)}.\label{pss6}
\end{align}
As an illustration, note that if $r_1^{(2,3)}, r_2^{(1,3)}$ are non-zero elements of $\Bbbk$, then we obtain $\beta=\gamma$ and $\gamma^{-1}  = \alpha^{-1} $, from (\ref{cattt}) and  (\ref{doggg}), respectively. So, (\ref{delfi}) implies that $r_3^{(1,3)} = r_3^{(1,2)}$. Of course, if $r_1^{(2,3)}, r_2^{(1,3)}=0$ it is not necessarily true that $\gamma^{-1}  = \alpha^{-1} $ and $\beta = \alpha^{-1}$.
\begin{examples}\label{chapeco}
\begin{enumerate}
\item \textit{Woronowicz algebra $\cW_{\nu}(\mathfrak{sl}(2,\Bbbk))$.} This $\Bbbk$-algebra was
introduced by Woronowicz in \cite{Woronowicz1987}. It is generated
by the indeterminates $x,y,z$ subject to the relations
$xz-\nu^4zx=(1+\nu^2)x,\ xy-\nu^2yx=\nu z,\
zy-\nu^4yz=(1+\nu^2)y$, where $\nu \in\ \Bbbk\ \backslash\ \{0\}$ is not a root of unity. Under certain
conditions on $\nu$ (see at the end of the computations), this algebra is a 3-dimensional skew polynomial algebra. Let us see the details. Let $x_1:=x, x_2:=y$
and $x_3:=z$. We have the relations $x_3x_1=\nu^{-4}x_1x_3-\nu^{-4}(1+\nu^2)x_1,\
x_2x_1=\nu^{-2}x_1x_2-\nu^{-1}x_3,\ x_3x_2=\nu^4x_2x_3+(1+\nu^{2})x_2$. If $x_1\prec x_2\prec x_3$, then $(Q,\preceq_{\rm deglex})$ is a
skew reduction system and
\begin{align*}
{\rm stred}_Q&(x_3f_{21})=x_3(\nu^{-2}x_1x_2-\nu^{-1}x_3)\\
&=\nu^{-2}x_3x_1x_2-\nu^{-1}x_3^2\\
&=\nu^{-2}(\nu^{-4}x_1x_3-\nu^{-4}(1+\nu^2)x_1)x_2-\nu^{-1}x_3^2\\
&=\nu^{-6}x_1x_3x_2-\nu^{-6}(1+\nu^2)x_1x_2-\nu^{-1}x_3^2\\
&=\nu^{-6}x_1(\nu^4x_2x_3+(1+\nu^2)x_2)-\nu^{-6}(1+\nu^2)x_1x_2-\nu^{-1}x_3^2\\
&= \nu^{-2}x_1x_2x_3+\nu^{-6}(1+\nu^2)x_1x_2-\nu^{-6}(1+\nu^2)x_1x_2-\nu^{-1}x_3^2\\
&= \nu^{-2}x_1x_2x_3-\nu^{-1}x_3^2,
\end{align*}
while,
\begin{align*}
{\rm stred}_Q(f_{32}x_1)&=(\nu^4x_2x_3+(1+\nu^2)x_2)x_1\\
&=\nu^4x_2x_3x_1+(1+\nu^2)x_2x_1\\
&=\nu^4x_2(\nu^{-4}x_1x_3-\nu^{-4}(1+\nu^2)x_1)+(1+\nu^2)(\nu^{-2}x_1x_2-\nu^{-1}x_3)\\
&= x_2x_1x_3-(1+\nu^2)x_2x_1+(1+\nu^2)\nu^{-2}x_1x_2-(1+\nu^2)\nu^{-1}x_3\\
&=(\nu^{-2}x_1x_2-\nu^{-1}x_3)x_3-(1+\nu^2)(\nu^{-2}x_1x_2-\nu^{-1}x_3)\\
& + (1+\nu^2)\nu^{-2}x_1x_2-\nu^{-1}(1+\nu^2)x_3\\
&=\nu^{-2}x_1x_2x_3-\nu^{-1}x_3^2-\nu^{-2}(1+\nu^2)x_1x_2 + (1+\nu^2)\nu^{-1}x_3\\
& + (1+\nu^2)\nu^{-2}x_1x_2-\nu^{-1}(1+\nu^2)x_3\\
&=\nu^{-2}x_1x_2x_3-\nu^{-1}x_3^2.
\end{align*}
Thus, Theorem \ref{GomezTorrecillas2Theorem 4.7}
implies that $\cW_{\nu}(\mathfrak{sl}(2,\Bbbk))$ is a 3-di\-men\-sio\-nal skew polynomial algebra for any value of $\nu \in \Bbbk\ \backslash\ \{0\}$.
\item \textit{Dispin algebra $\cU(osp(1,2))$.} This $\Bbbk$-algebra is generated by the variables $x,y,z$ subjected to the relations $yz-zy=z,\ zx+xz=y,\ xy-yx=x$. Let $x_1:=x, x_2:=y$ and $x_3:=z$. We consider $x_1\prec x_2\prec x_3$. Then
$(Q,\preceq_{\rm deglex})$ is a skew reduction system. Relations
defining this algebra are $x_3x_2=x_2x_3-x_3,\ x_3x_1=-x_1x_3+x_2$, and $
x_2x_1=x_1x_2-x_1$. Following Theorem \ref{GomezTorrecillas2Theorem 4.7}, we have
\begin{align*}
{\rm stred}_Q(x_3f_{21})&=x_3(x_1x_2-x_1)\\
&= x_3x_1x_2-x_3x_1\\
&=(-x_1x_3+x_2)x_2-(-x_1x_3+x_2)\\
&=-x_1x_3x_2+x_2^2+x_1x_3-x_2\\
&=-x_1(x_2x_3-x_3)+x_2^2+x_1x_3-x_2\\
&=-x_1x_2x_3+x_1x_3+x_2^2+x_1x_3-x_2,
\end{align*}
and
\begin{align*}
{\rm stred}_Q(f_{32}x_1)&=(x_2x_3-x_3)x_1\\
&=x_2x_3x_1-x_3x_1\\
&=x_2(-x_1x_3+x_2)-(-x_1x_3+x_2)\\
&=-x_2x_1x_3+x_2^2+x_1x_3-x_2\\
&=-(x_1x_2-x_1)x_3+x_2^2+x_1x_3-x_2\\
&=-x_1x_2x_3+x_1x_3+x_2^2+x_1x_3-x_2.
\end{align*}
We can see that ${\rm stred}_Q(x_3f_{21})={\rm
stred}_Q(f_{32}x_1)$, so Theorem \ref{GomezTorrecillas2Theorem 4.7}
guarantees that $\cU(osp(1,2))$ is a 3-dimensional skew polynomial algebra.
\item If we consider the $\Bbbk$-algebra $A$ generated by the variables $x, y, x$ subjected to the relations $yx = \alpha xy + x$, $zx = \beta xz + z$, $zy=yz$, with $\beta, \alpha \in \Bbbk\ \backslash\ \{0\}$, then the set of variables $\{x, y, x\}$ is not a PBW basis for $A$. Consider the identification $x_1:=x,\ x_2:= y, x_3:= z$. Then, the algebra $A$ is expressed by the relations $x_2x_1 = \alpha x_1x_2 + x_1,\ x_3x_1 = \beta x_1x_3 + x_3,\ x_3x_2 = x_2x_3$, and hence
\begin{align*}
{\rm stred}_Q(x_3f_{21}) = &\ x_3(\alpha x_1x_2 + x_1) = \alpha x_3x_1x_2 + x_3x_1\\
= &\ \alpha (\beta x_1x_3 + x_3)x_2 + \beta x_1x_3 + x_3\\
= &\ \beta \alpha x_1x_3x_2 + \alpha x_3x_2 + \beta x_1x_3 + x_3\\
= &\ \beta \alpha x_1x_2x_3 + \alpha x_2x_3 + \beta x_1x_3 + x_3\\
{\rm stred}_Q(f_{32}x_1) = &\ x_2x_3x_1\\
= &\ x_2(\beta x_1x_3 + x_3) = \beta x_2x_1x_3 + x_2x_3\\
= &\ \beta(\alpha x_1x_2 + x_1)x_3 + x_2x_3\\
= &\ \beta \alpha x_1x_2x_3 + \beta x_1x_3 + x_2x_3.
\end{align*}
Since ${\rm stred}_Q(x_3f_{21}) \neq {\rm stred}_Q(f_{32}x_1)$, Theorem \ref{GomezTorrecillas2Theorem 4.7} guarantees that the set $\{x, y, z\}$ is not a PBW basis for the algebra $A$. Considering the notation in Remark \ref{Gererere}, we observe that $r_0^{(1,2)} = r_2^{(1,2)} = r_3^{(1,2)} = r_0^{(1,3)} = r_1^{(1,3)} = r_2^{(1,3)} = r_0^{(2,3)} = r_1^{(2,3)} = r_2^{(2,3)} = r_3^{(2,3)} = 0$, and $r_1^{(1,2)} = r_3^{(1,3)} = \alpha^{-1} = 1$. In particular, expression (\ref{delfi}) impose that $1=0$, which of course is false. This illustrates why the set $\{x, y, z\}$ is not a PBW basis over $\Bbbk$ for the algebra $A$.
\end{enumerate}
\end{examples}
\subsection*{Acknowledgment}
The first author is supported by Grant HERMES CODE 30366, Departamento de Matem\'aticas, Universidad Nacional de Colombia, Bogot\'a.


\end{document}